\documentclass[a4paper,12pt]{article}

\usepackage[latin1]{inputenc}
\usepackage{indentfirst}	
\usepackage{fancyhdr}
\usepackage[english]{babel} 
\usepackage{authblk} 
\usepackage{url} 
\usepackage[parfill]{parskip}
\usepackage{graphicx}
\usepackage{amssymb}
\usepackage{amsmath} 
\usepackage{latexsym}
\usepackage{amsthm}
\usepackage{mathrsfs} 

\begingroup
\makeatletter
\@for\theoremstyle:=definition,remark,plain\do{%
	\expandafter\g@addto@macro\csname th@\theoremstyle\endcsname{%
		\addtolength\thm@preskip\parskip
	}%
}
\endgroup

\usepackage[a4paper, margin=3cm]{geometry}

\allowdisplaybreaks

\usepackage[final,color]{showkeys}	
\usepackage{seqsplit}
\usepackage{xstring}
\renewcommand{\showkeyslabelformat}[1]{
	\noexpandarg
	\StrSubstitute{#1}{ }{~}[\TEMP]
	\fbox{\parbox{0.1\linewidth}{\normalfont\small\ttfamily\expandafter\seqsplit{\expandafter\TEMP}} 
	}}        
\definecolor{refkey}{rgb}{0.6,0.7,0.6}
\definecolor{labelkey}{rgb}{0.6, 0.7, 0.6}

\usepackage{tikz-cd} 
\usetikzlibrary{matrix,arrows}

\theoremstyle{plain} 
\newtheorem{thm}{Theorem}[section]

\newtheorem{lem}[thm]{Lemma}
\newtheorem{prop}[thm]{Proposition}
\newtheorem{claim}[thm]{Claim}

\theoremstyle{definition} 
\newtheorem{defn}[thm]{Definition}

\newtheorem{defn-prop}[thm]{Definition-Proposition}
\newtheorem{rem}[thm]{Remark}
\newcommand{\val}{\operatorname{val}}

\begin{document}

\thispagestyle{empty}
\begin{center}
	
\Large\textbf{The isotopy problem for the phase tropical line}

\normalsize\text{Raffaele Caputo}

\normalsize\text{University of Hamburg}

\end{center}

\begin{abstract}
 Let $H$ be the complex line $1+z_{1}+z_{2}=0$ in $(\mathbb{C}^{*})^{2}$ and $H_{trop}$ the associated phase tropical line. 
 We show that $H$ and $H_{trop}$ are isotopic as topological submanifolds, by explicitly constructing an isotopy map.
\end{abstract}

\section{Introduction}

Given an $n$-dimensional complex hyperplane $H\subset(\mathbb{C}^{*})^{n}$, one can construct the associated phase tropical hyperplane $H_{trop}\subset(\mathbb{C}^{*})^{n}$. This is a cell complex whose projection to $\mathbb{R}^{n}$ is the tropical hyperplane.
In $\cite{Kerr}$ and $\cite{Kim}$ it has been shown that in any dimension  these two objects are homeomorphic. In the present work we show that for $n=2$ one can get a stronger result, namely one can continuously deform $H$ into $H_{trop}$ via homeomorphic objects in $(\mathbb{C}^{*})^{2}$. A generalization of this result in higher dimensions can be found in the upcoming work $\cite{RuZh}$ by H.~Ruddat and I.~Zharkov. We prove the following:
\begin{thm}\label{thm:1.1}
Let $H\subset(\mathbb{C}^{*})^{2}$ the complex line $1+z_{1}+z_{2}=0$ and $H_{trop}\subset(\mathbb{C}^{*})^{2}$ its associated phase tropical line. Let $i_{H}:H\hookrightarrow(\mathbb{C}^{*})^{2}$ be the canonical embedding.
There exists a continuous map: 
\[\Psi:H\times [0,1]\longrightarrow (\mathbb{C}^{*})^{2}\] 
such that, for $t\in [0,1]$, the family of maps:  
\[\Psi_{t}:H\longrightarrow (\mathbb{C}^{*})^{2}, P\mapsto \Psi(P,t),\]
has the following properties: 
\begin{enumerate}
	\item $\Psi_{0}=i_{H}$;
  \item $\Psi_{1}(H)=H_{trop}$;
	\item $\Psi_{t}$ is a homeomorphism onto the image, for each $t\in\left[0,1\right]$.
\end{enumerate}
\end{thm} 
Throughout, we identify $(\mathbb{C}^{*})^{2}$ with $\mathbb{R}^{2}\times (S^{1})^{2}$.
We start in Section $\ref{sec: phase tropical}$ by defining the phase tropical line following $\cite{Mik}$ and $\cite{Nis}$. We give a practical description of $H_{trop}$ as fibred over its image under the canonical projection $\pi:\mathbb{R}^{2}\times (S^{1})^{2}\rightarrow \mathbb{R}^{2}$. In Section $\ref{sec: complex line}$ we do the same thing for the complex line.   
In Section $\ref{sec: preparation}$ we prepare $H$ and $H_{trop}$ for the construction of the isotopy map. To do that, we cut both $H$ and $H_{trop}$ in three subsets, which we call $H_i$ and $H_{itrop}$, for $i\in \{1,2,3\}$, respectively (see Figure~5, Figure 6, Figure 7 and Definition $\ref{defn:44}$).  In Lemma $\ref{lem:4}$ we show that, both in $H$ and $H_{trop}$, these three subsets are swapped in a cyclic manner by the following automorphism of order $3$ of $\mathbb{R}^{2}\times (S^{1})^{2}$:
\begin{equation}
\begin{split}
\lambda:\mathbb{R}^{2}\times (S^{1})^{2}&\longrightarrow \mathbb{R}^{2}\times (S^{1})^{2}\\
(x,y,\varphi,\psi)&\longmapsto (-y,x-y,-\psi+2\pi,\varphi-\psi+2\pi)\\
\end{split}.
\end{equation}
More precisely, we have $\lambda(H_i)=H_j$ and $\lambda(H_{itrop})=H_{jtrop}$, with $j=i+1$, for $i\in\{1,2\}$, and $j=i\mod{2}$, for $i=3$. Therefore, we focus our attention on the subsets $H_{1}\subset H$ and $H_{1trop}\subset H_{trop}$ in order to construct an isotopy between the canonical embedding of $H_{1}$ in $(\mathbb{C}^{*})^{2}$ and an embedding of $H_{1}$ in $(\mathbb{C}^{*})^{2}$ whose image is $H_{1trop}$. We call this isotopy map $\Phi^{t}_{1}$. To do that we subdivide $H_{1}$ and $H_{1trop}$ in two parts (see Figure 8 and Definition $\ref{defn:45}$).
Finally, in Section $\ref{sec: isotopy}$, we explicitly write the isotopy map $\Phi^{t}_{1}$ and then we prove Theorem $\ref{thm:1.1}$.

\text{\textit{Acknowledgements}.} I would like to thank Bernd Siebert for inspiring discussions and insights on the subject.

\section{The phase tropical line}
\label{sec: phase tropical}

Let $K$ be the field of the (real-power) Puiseux series:
\[
K:=\left\{\sum_{j\in I}c_{j}t^{j}|c_{j}\in\mathbb{C}, I\subset\mathbb{R} \text{ well ordered }\right\}.
\]
Set $K^{*}:=K\setminus\left\{0\right\}$. Since $I$ is a well ordered set, any element $z\in K$ can be written as $z=c_{j_0}t^{j_0}+\text{higher terms}$. Set $r_0:=|c_{j_0}|$ and $\varphi_0:=\arg c_{j_0}$. The field $K$ comes with a non-Archimedean valuation. Recall that a non-Archimedian valuation
is a function $\val : K^{*} \rightarrow \mathbb{R}$, such that $\val(a+b) \leq\max(\val(a), \val(b))$ 
and $\val(ab) = \val(a) +  \val(b)$. In our case, the valuation is
defined by $\val(r_0 e^{i\varphi_0}t^{j_{0}}+ \text{ higher terms}):=-j_0$. We notice that in the definition of $K$ we use irrational as well as rational powers in the Puiseux series in order to make the valuation surjective. If we keep the phase $\varphi_0$ of $c_{j_0}$ we can lift the valuation map to $\mathbb
{C}^{*}$. Namely, we can define a map
\begin{equation}
\begin{split}
 \val_{\mathbb{C}}:K^{*}&\longrightarrow \mathbb{C}^{*}\\
 r_0 e^{i\varphi_0}t^{j_{0}}+ \text{ higher terms}&\longmapsto e^{-j_{0}+i\varphi_0}\\ 
\end{split}.
\end{equation}

This map extends componentwise to a map $\val^{2}_{\mathbb{C}}:(K^{*})^{2}\longrightarrow (\mathbb{C}^{*})^{2}$.
We are ready to define the phase tropical line following $\cite{Mik}$ and $\cite{Nis}$.
Let $H$ be the complex line in $(\mathbb{C}^{*})^{2}$:
\[H:=\left\{(z_{1},z_{2})\in(\mathbb{C}^{*})^{2}|1+z_{1}+z_{2}=0\right\}.\]
Let $f:=1+z_{1}+z_{2}\in K[z_{1},z_{2}]$ and consider its zero set $Z(f)$ in  $(K^{*})^{2}$. 
\begin{defn}\textit{The phase tropical line} $H_{trop}$ associated to $H$ is 
\[
H_{trop}:=\overline{\val^{2}_{\mathbb{C}}(Z(f))}\subset (\mathbb{C}^{*})^{2}.
\]
\end{defn}

We want to find a practical way to describe and understand $H_{trop}$. To do that, we first notice that if $f:=1+z_{1}+z_{2}\in K[z_{1},z_{2}]$, then its zero set $Z(f)\subset (K^{*})^{2}$ can be written as the set of all pairs in $(K^{*})^{2}$ of the form $(z,-1-z)$ for  $z\in  K^{*}\setminus \left\{-1\right\}$. Thus, if we define the map:
\[\nu:K^{*}\setminus \left\{-1\right\}\longmapsto (\mathbb{C}^{*})^{2}, z\longmapsto(\val(z),\val(-1-z)),\] 
then 
\[H_{trop}=\overline{\nu(K^{*}\setminus \left\{-1\right\})}.\]
Moreover, we identify $(\mathbb{C}^{*})^{2}$ with $\mathbb{R}^{2}\times (S^{1})^{2}$ via the homeomorphism:
\[h:(\mathbb{C}^{*})^{2}\longrightarrow \mathbb{R}^{2}\times (S^{1})^{2}, (z_{1},z_{2})\longmapsto (\ln|z_{1}|,\ln|z_{2}|,\arg z_{1},\arg z_{2}).\]
Composing $\val_{\mathbb{C}}$ with $h$ we get a map:
\[\val_{\mathbb{C}}':K^{*}\stackrel{\val_{\mathbb{C}}}{\longrightarrow}\mathbb{C}^{*}\stackrel{h}{\longrightarrow}\mathbb{R}\times S^{1}\] \[r_0e^{i\varphi_0} t^{j_{0}}+ \text{ higher terms}\longmapsto(-j_{0}, \varphi_0).\]
Analogously, composing $\nu$ with $h\times h$ we get a map
\[\nu':K^{*}\setminus \left\{-1\right\}\stackrel{\nu}{\longrightarrow}(\mathbb{C}^{*})^{2} \stackrel{h\times h}{\longrightarrow}\mathbb{R}^{2}\times (S^{1})^{2}\] 
\[z\longmapsto(\val_{\mathbb{C}}'(z), \val_{\mathbb{C}}'(-1-z)).\]
 Thus, if we embed $H_{trop}$ in $\mathbb{R}^{2}\times (S^{1})^{2}$, we can write \[H_{trop}=\overline{\nu'(K^{*}\setminus \left\{-1\right\})}.\]
  We see that to describe the phase tropical line we need to analyze the output of $\val_{\mathbb{C}}'~(-1-z)$, for $z\in K^{*}\setminus \left\{-1\right\}$. That is, we need to determine the term with the smallest exponent in the sum $-1-z$. To do that, we subdivide $K^{*}\setminus \left\{-1\right\}$ in four subsets. Namely, writing $z=r_0e^{i\varphi_0} t^{j_{0}}+ \text{ higher terms}$, for $z\in K^{*}\setminus \left\{-1\right\}$, we set $P_{1}:=\left\{z\in K^{*}\setminus \left\{-1\right\}|j_{0}> 0\right\}$, $P_{2}~:=\left\{z\in K^{*}\setminus \left\{-1\right\}|j_{0}=0, r_{0}e^{i\varphi_0}= -1\right\}$, $P_{3}:=\left\{z\in K^{*}\setminus \left\{-1\right\}| j_{0}< 0\right\}$ and $P_{4}:=\left\{z\in K^{*}\setminus \left\{-1\right\}| j_{0}= 0, r_{0}e^{i\varphi_0}\neq -1\right\}$. Now, if $z\in P_1$, then in the sum $-1-z$ the term with the smallest exponent is $-1$, thus $\val_{\mathbb{C}}'~(-1-z)=-1$. Hence,
\[\nu'(P_{1})=\left\{(x,0,\varphi,\pi)\in\mathbb{R}^{2}\times (S^{1})^{2}|x\in\mathbb{R}_{<0}, \varphi\in S^{1}\right\}.\]
If $z\in P_2$, then in the sum $-1-z$ the term with the smallest exponent is $\alpha t^{j}$, for any $j>0$ and $\alpha\in\mathbb{C}^{*}$. Thus, 
\[\nu'(P_{2})=\left\{(0,y,\pi,\psi)\in\mathbb{R}^{2}\times (S^{1})^{2}|y\in\mathbb{R}_{<0}, \psi\in S^{1}\right\}.\]
Analogously, if $z\in P_3$, then in the sum $-1-z$ the term with the smallest exponent is $-r_{0}e^{i\varphi_0}t^{j_{0}}$. Thus, 
\[\nu'(P_{3})=\left\{(x,x,\varphi,\varphi+\pi)\in\mathbb{R}^{2}\times (S^{1})^{2}|x\in\mathbb{R}_{>0}, \varphi\in S^{1}\right\}.\]
Finally,  if $z\in P_4$, then in the sum $-1-z$ the term with the smallest exponent is $-1-r_{0}e^{i\varphi_0}t^{j_{0}}$. Thus, 
\[\nu'(P_{4})=\left\{(0,0,\varphi,\arg(-1-re^{i\varphi}))\in\mathbb{R}^{2}\times (S^{1})^{2}|r\in\mathbb{R}_{>0}, \varphi\in S^{1}\right\}.\]
Taking the closure in $\mathbb{R}^{2}\times (S^{1})^{2}$ of the union of the subsets $\nu'(P_{i})$, for $i=1,2,3,4$, we get $H_{trop}$.

\begin{defn}
Let $H$ be the complex line and $H_{trop}$ the phase tropical line. Let $\pi:\mathbb{R}^{2}\times (S^{1})^{2}\rightarrow \mathbb{R}^{2}$ and $Arg:\mathbb{R}^{2}\times (S^{1})^{2}\rightarrow (S^{1})^{2}$ be the projections. Then $\pi(H)$ and $\pi(H_{trop})$ are called respectively the \textit{amoeba} and the \textit{tropical amoeba} of $H$. 
Moreover, $Arg(H)$ and $Arg(H_{trop})$ are called respectively the \textit{coamoeba} and \textit{tropical coamoeba} of $H$. 
\end{defn}

The following proposition can be found in a more general form in $\cite{Nis}$.
\begin{prop}\label{prop:30}
Let $H$ be the complex line and $H_{trop}$ the phase tropical line, then
the tropical coamoeba equals the closure of the coamoeba.
\end{prop}
The tropical amoeba of $H$ is given by the union of three half lines:
\[\pi(H)=\left\{x\leq 0,y=0\right\}\cup\left\{x=0,y\leq 0\right\}\cup\left\{x=y,y\geq0\right\}.\]
\begin{center}
\includegraphics[scale=0.5]{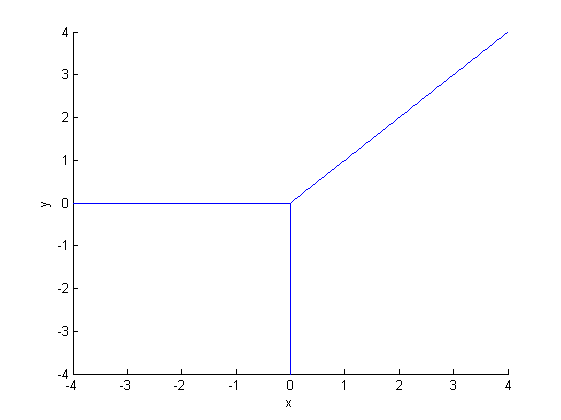}
\begin{center}
\textbf{Figure 1}: The tropical amoeba.
\end{center}
\end{center}
The tropical coamoeba, seen in the fundamental domain $\left[0,2\pi\right]^{2}$, can be pictured as follows:
\begin{center}
\includegraphics[scale=0.5]{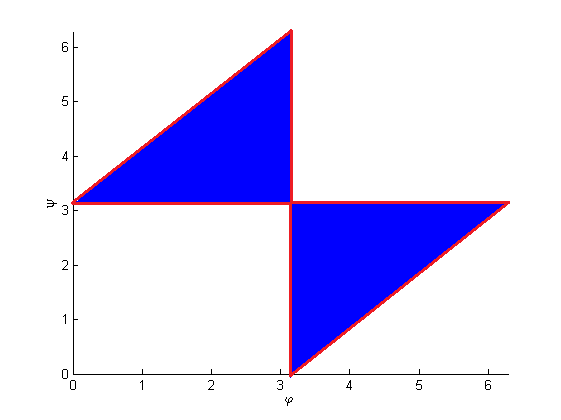}
\begin{center}
\textbf{Figure 2}: The tropical coamoeba.\end{center}
\end{center}
The picture is justified by Proposition $\ref{prop:30}$ and Proposition $\ref{prop:3}$.
\begin{rem}
	We notice that in the definition of the phase tropical line we need to take the closure, for without it the phase
	tropical line would omit the six line segments that are in the closure of the coamoeba but not in the coamoeba (see Figure 4).
\end{rem}
 If we consider the projection $\pi|_{H_{trop}}:H_{trop}\longrightarrow \pi(H_{trop})$, from the above description of the phase tropical line we see that the fibre of $\pi|_{H_{trop}}$ is given by $\left\{(\varphi,\pi)|\varphi\in S^{1}\right\}$ over $x< 0$, by $\left\{(\pi,\psi)|\psi\in S^{1}\right\}$ over $y< 0$, by $\left\{(\varphi,\varphi+\pi)|\varphi\in S^{1}\right\}$ over $x=y, y>0$ and by the whole tropical coamoeba $Arg(H_{trop})$ over the point $(0,0)$.

\section{The complex line}
\label{sec: complex line}
 
We give a detailed description of the amoeba and coamoeba of the complex line $H$. In Proposition $\ref{prop:2}$ we describe $H$ by describing the fibres of $\pi|_{H}:H\longrightarrow \pi(H)$. 
\begin{lem}\label{lem:1}
Let $H$ be the complex line:
\[H=\left\{(x,y,\varphi,\psi)\in\mathbb{R}^{2}\times(S^{1})^{2}|e^{x+i\varphi}+e^{y+i\psi}+1=0 \right\}.\]
If $(x,y,\varphi,\psi)\in H$, then:
\begin{equation}\label{eq:1}
\begin{cases}
e^{2x}=1+2e^{y}\cos\psi+e^{2y}\\e^{2y}=1+2e^{x}\cos\varphi+e^{2x}
\end{cases}.
\end{equation}
Conversely, if $(x,y,\varphi,\psi)\in\mathbb{R}^{2}\times(S^{1})^{2}$ satisfies $\eqref{eq:1}$, then: 
\[\left\{(x,y,\varphi,\psi),(x,y,2\pi-\varphi,2\pi-\psi)\right\}\text{ or }\left\{(x,y,2\pi-\varphi,\psi),(x,y,\varphi,2\pi-\psi)\right\}\]
is a subset of $H$.
\end{lem}

\proof
From the defining equation of $H$ we get:
\[
\begin{cases}
e^{x}=|-1-e^{y+i\psi}|\\e^{y}=|-1-e^{x+i\varphi}|
\end{cases}
\]
which, elevating both sides of each equation to the square power and using the trigonometric expression for the complex exponential function, can be seen to be equivalent to:
\begin{equation}\label{eq:2}
\begin{cases}
(1+e^{y}\cos\psi)^{2}+(e^{y}\sin\psi)^{2}=e^{2x}\\(1+e^{x}\cos\varphi)^{2}+(e^{x}\sin\varphi)^{2}=e^{2y}
\end{cases}
\end{equation}
which, in turn, using trigonometric identities, can be seen to be equivalent to $\eqref{eq:1}$.

Conversely, let $(x,y,\varphi,\psi)\in\mathbb{R}^{2}\times (S^{1})^{2}$ satisfy $\eqref{eq:1}$.
Then, we immediately see that $(x,y,2\pi-\varphi,2\pi-\psi)$, $(x,y,2\pi-\varphi,\psi)$ and $(x,y,\varphi,2\pi-\psi)$ are also solutions of $\eqref{eq:1}$.
Now, we show that only one of the two sets:
\[\left\{(x,y,\varphi,\psi),(x,y,2\pi-\varphi,2\pi-\psi)\right\}\text{ , }\left\{(x,y,2\pi-\varphi,\psi),(x,y,\varphi,2\pi-\psi)\right\},\]
both if they coincide, also satisfies $e^{x+i\varphi}+e^{y+i\psi}+1=0$.
To see this, we notice that $\eqref{eq:1}$ implies $e^{2x}+e^{2y}=2+2e^{y}\cos\psi+2e^{x}\cos\varphi+e^{2y}+e^{2x}$, which in turn implies:
\begin{equation}\label{eq:3}
e^{x}\cos\varphi+ e^{y}\cos\psi+1=0. 
\end{equation}
On the other hand $\eqref{eq:1}$ implies $\eqref{eq:2}$, and hence:
\begin{equation}\label{eq:4}
(1+e^{y}\cos\psi)^{2}+(e^{y}\sin\psi)^{2}=e^{2x}.
\end{equation}

Substituting $\eqref{eq:3}$ in $\eqref{eq:4}$, we get $(e^{y}\sin\psi)^{2}=e^{2x}-e^{2x}\cos^{2}\varphi$, that is:
\begin{equation}\label{eq:5}
(e^{y}\sin\psi)^{2}=(e^{x}\sin\varphi)^2.
\end{equation}
We see that only one of the two above sets, both if they coincide, satisfies:
\begin{equation}\label{eq:5a}
e^{y}\sin\psi+e^{x}\sin\varphi=0,
\end{equation}
while the other satisfies:
\begin{equation}
e^{y}\sin\psi-e^{x}\sin\varphi=0.
\end{equation}
Thus, the subset satisfying $\eqref{eq:5a}$ also satisfies:
\begin{equation} 
\begin{cases}
e^{x}\cos\varphi+ e^{y}\cos\psi+1=0 
\\ e^{x}\sin\varphi+ e^{y}\sin\psi=0 
\end{cases} ,
\end{equation}
which is equivalent to $e^{x+i\varphi}+e^{y+i\psi}+1=0$.

\endproof

\begin{prop}
Let $H$ be the complex line and $\pi(H)\subset \mathbb{R}^{2}$ its amoeba. Let ${(x,y)\in\mathbb{R}^{2}}$. Then, $(x,y)\in\pi(H)$ if and only if:
\begin{equation}\label{eq:6}
\begin{cases}
 e^{x}-e^{y}\leq 1
\\ e^{y}-e^{x}\leq 1
\\ e^{x}+e^{y}\geq 1
\end{cases}.
\end{equation}
\end{prop}
\proof
Let $(x,y)\in\mathbb{R}^{2}$. By Lemma $\ref{lem:1}$, $(x,y)\in\pi(H)$ if and only if we can find $(\varphi,\psi)\in (S^{1})^{2}$ for which $\eqref{eq:1}$ holds. Now, once we fix $(x,y)\in\mathbb{R}^{2}$, $\eqref{eq:1}$ has a solution in $(S^{1})^{2}$ if and only if: 
\[
\begin{cases}
-1\leq \frac{e^{2x}-e^{2y}-1}{2e^{y}}\leq 1\\-1\leq\frac{e^{2y}-e^{2xy}-1}{2e^{x}}\leq 1
\end{cases},
\]
which is equivalent to:

\begin{equation}\label{eq:7}
\begin{cases}
(e^{x}-e^{y}-1)(e^{x}+e^{y}+1)\leq 0\\(e^{x}-e^{y}+1)(e^{x}+e^{y}-1)\geq 0\\(e^{y}-e^{x}-1)(e^{y}+e^{x}+1)\leq 0\\(e^{y}-e^{x}+1)(e^{y}+e^{x}-1)\geq 0
\end{cases}.
\end{equation}
We easily see that $\eqref{eq:7}$ holds if and only if $\eqref{eq:6}$ holds.

\endproof

\begin{center}
\includegraphics[scale=0.6]{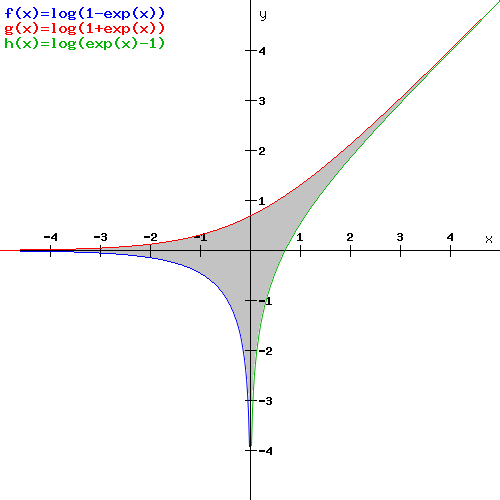}
\begin{center}
	\textbf{Figure 3}: The amoeba of $H$.
\end{center}
\end{center}

\begin{defn}
We denote by $H_{\circ}$ the subset of all points $(x,y,\varphi,\psi)\in H$ whose image under $\pi$ lies in the interior of the amoeba, that is $(x,y)$ satisfies:
\[
\begin{cases}
 e^{x}-e^{y}< 1
\\ e^{y}-e^{x}< 1
\\ e^{x}+e^{y}> 1
\end{cases}.
\] 
We denote by $H_{\partial }$ the subset of all points $(x,y,\varphi,\psi)\in H$ whose image under $\pi$ lies on the boundary of the amoeba, that is $(x,y)$ either satisfies $ e^{x}-e^{y}=1$, $e^{y}-e^{x}=1$ or $e^{x}+e^{y}=1$.
\end{defn}

Using Lemma $\ref{lem:1}$ one can prove the following:
\begin{prop}\label{prop:2}
Let $(x,y,\varphi,\psi)\in H$. Then:
\[
(\varphi,\psi)=
\begin{cases}
(\pi,0) &\text{ if and only if } e^{x}-e^{y}=1\\
(0,\pi) &\text{ if and only if } e^{y}-e^{x}=1\\
(\pi,\pi) &\text{ if and only if } e^{x}+e^{y}=1
\end{cases}.
\]
Let $(x,y)\in\pi(H_{\circ})$, then:
\[\pi|_{H}^{-1}(x,y)=\left\{(x,y,\varphi,\psi),(x,y,2\pi-\varphi,2\pi-\psi)\right\},\]
for some $(\varphi,\psi)\in (S^{1})^{2}\setminus \left\{(0,0),(0,\pi),(\pi,0),(\pi,\pi)\right\}$.
\end{prop}
\proof
Let $P:=(x,y,\varphi,\psi)\in H$, thus $P$ satisfies $1+e^{x+i\varphi}+e^{y+i\psi}=0$. Now, $1+e^{x+i\varphi}+e^{y+i\psi}=0$ is equivalent to $1-e^{x}+e^{y}=0$  if and only if $e^{i\varphi}=-1$ and $e^{i\psi}=1$, that is if and only if $\varphi=\pi$ and $\psi=0$. In the same way, we get that $P$ satisfies $e^{y}-e^{x}=1$ if and only if $\varphi=0$ and $\psi=\pi$. Finally, $P\in H$ satisfies $e^{x}+e^{y}=1$ if and only if $\varphi=\pi$ and $\psi=\pi$.
Now, let us assume $P\in H_{\circ}$. By Lemma $\ref{lem:1}$, we immediately get:

\[\pi|_{H}^{-1}(x,y)=\left\{(x,y,\varphi,\psi),(x,y,2\pi-\varphi,2\pi-\psi)\right\}.\]

By the first part of the proof, $(\varphi,\psi)\in (S^{1})^{2}\setminus \left\{(0,\pi),(\pi,0),(\pi,\pi)\right\}$. Moreover, $(\varphi,\psi)\neq (0,0)$ as $1+e^{x+i\varphi}+e^{y+i\psi}=0$ computed in $(\varphi,\psi)= (0,0)$ has no solution in $\mathbb{R}^{2}$.
\endproof

\begin{defn}
We define two subsets of $H$, namely \[H^{+}:=\left\{(x,y,\varphi,\psi)\in H|0\leq\varphi\leq \pi\right\} \text{ and }H^{-}:=\left\{(x,y,\varphi,\psi)\in H|\pi\leq\varphi\leq 2\pi\right\}.\]
\end{defn}

\begin{prop}\label{prop:33}
The maps $\pi_{+}:=\pi|_{H^{+}}$ and $\pi_{-}:=\pi|_{H^{-}}$ are homeomorphisms onto their images.
\end{prop}
\proof
By Proposition $\ref{prop:2}$ the maps are bijective. They are continuous being restrictions of a continuous map. Moreover, since $(S^{1})^{2}$ is compact, the projection $\pi$ is a closed map, therefore $\pi_{+}$ and $\pi_{-}$ are also closed maps being restrictions of $\pi$ to closed subsets.
\endproof

\begin{lem}\label{lem:2}
Let $(x,y,\varphi,\psi)\in H$. Then, $(x,y,\varphi,\psi)\in H_{\circ}$ if and only if
\begin{equation}\label{eq:8}
\begin{cases}
 e^{x}=-\frac{\sin\psi}{\sin(\psi-\varphi)}
\\ e^{y}=+\frac{\sin\varphi}{\sin(\psi-\varphi)}
\end{cases} 
\end{equation}
\end{lem}
\proof
Let $(x,y,\varphi,\psi)\in H$. By Proposition $\ref{prop:2}$, $(x,y,\varphi,\psi)\in H_{\circ}$ if and only if  $(\varphi,\psi)\notin\left\{(0,0),(0,\pi),(\pi,0),(\pi,\pi)\right\}$. The defining equation of $H$ is equivalent to the following system of equations:
\begin{equation}\label{eq:9}
\begin{cases}
 e^{x}\cos\varphi+ e^{y}\cos\psi+1=0 
\\ e^{x}\sin\varphi+ e^{y}\sin\psi=0 
\end{cases} .
\end{equation}
From $(\ref{eq:9})$ we see that if $\sin\psi\neq 0$, then $\sin\varphi\neq 0$. Now, 
\[
\begin{cases}
\sin\psi\neq 0\\
\sin\varphi\neq 0
\end{cases}
\]
if and only if $(\varphi,\psi)\notin\left\{(0,0),(0,\pi),(\pi,0),(\pi,\pi)\right\}$. Under this assumption, $\eqref{eq:9}$ is equivalent to:
\begin{equation} \label{eq:10}
\begin{cases}
 e^{x}\sin(\psi-\varphi)+\sin\psi=0
\\ e^{x}\sin\varphi+e^{y}\sin\psi=0
\end{cases} .
\end{equation}
From the first equation in $\eqref{eq:10}$, we see that if $\sin\psi\neq 0$, then also $\sin(\psi-\varphi)\neq 0$. Therefore, $\eqref{eq:10}$ is equivalent to $\eqref{eq:8}$.

\endproof

\begin{defn}
We define two subsets of $(S^{1})^{2}$, namely:
\[T_{1}:=\left\{(\varphi,\psi)\in(S^{1})^{2}|0<\varphi<\pi,\pi<\psi<\varphi+\pi\right\}\]
and
\[T_{2}:=\left\{(\varphi,\psi)\in(S^{1})^{2}|\pi<\varphi<2\pi,\varphi-\pi<\psi<\pi\right\}.\]
See Figure 4 below.
\end{defn}

\begin{prop}\label{prop:3}
The coamoeba of $H$ is:
\[Arg(H)=T_{1}\cup T_{2} \cup \left\{(0,\pi),(\pi,0),(\pi,\pi)\right\}.\]
Moreover, the map $Arg|_{H_{\circ}}:H_{\circ}\longrightarrow T_{1}\cup T_{2} $ is a homeomorphism.
\end{prop}
\proof
By Proposition $\ref{prop:2}$, if $(x,y,\varphi,\psi)\in H_{\partial}$, then its image under  $Arg$ is $(0,\pi)$,$(\pi,0)$ or $(\pi,\pi)$.
On the other hand, by Lemma $\ref{lem:2}$ if $(x,y,\varphi,\psi)\in H_{\circ}$, then $(x,y,\varphi,\psi)$ satisfies $\eqref{eq:8}$, hence $(\varphi,\psi)$ satisfies:
\begin{equation}
\begin{cases}
-\frac{\sin\psi}{\sin(\psi-\varphi)}>0
\\ +\frac{\sin\varphi}{\sin(\psi-\varphi)}>0.
\end{cases} 
\end{equation}
Therefore, if $\sin(\psi-\varphi)>0$, then $\sin\psi<0$ and $\sin\varphi>0$. In this case, $\sin\psi<0$ implies $\pi<\psi<2\pi$, and $\sin\varphi>0$ implies $0<\varphi<\pi$. Moreover, $\sin(\psi-\varphi)>0$ implies $0<\psi-\varphi<\pi$. Therefore, we get that $(\psi,\varphi)$ lies in $T_{1}$.
Conversely, if $(\varphi,\psi)$ lies in $T_{1}$, that is $0<\varphi<\pi$,$\pi<\psi<\varphi+\pi$, then we have $\sin\varphi>~0, \sin\psi<~0,\sin(\psi-\varphi)>0$, hence $\eqref{eq:8}$ has a $($unique$)$ solution.
The symmetry of $H$, $x\leftrightarrow y, \varphi\leftrightarrow\psi$, completes the proof of the first statement. 
In particular, we have shown that $Arg(H_{\circ})=T_{1}\cup T_{2}$. Now, $Arg|_{H_{\circ}}$ is continuous being the restriction of a continuous map. Its inverse map is:
\[Arg|_{H_{\circ}}^{-1}:T_{1}\cup T_{2}\longrightarrow H_{\circ}\]
\[(\varphi,\psi)\longmapsto (\ln(-\frac{\sin\psi}{\sin(\psi-\varphi)}),\ln(\frac{\sin\varphi}{\sin(\psi-\varphi)}),\varphi,\psi),\]
which is continuous as well.
\endproof

\begin{center}
\includegraphics[scale=0.6]{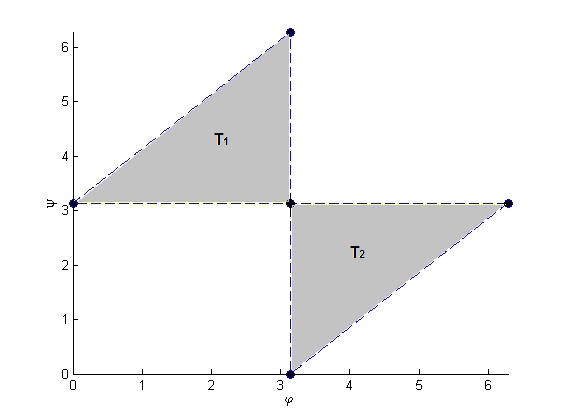}

\textbf{Figure 4}: The coamoeba of $H$.
\end{center}

\section{The isotopy preparation procedure}
\label{sec: preparation}

\begin{defn}\label{defn:4}
We subdivide $H$ in three parts. Namely:

\[H_{1}:=\left\{(x,y,\varphi,\psi)\in H|x\leq y, x\leq 0\right\}\]
\[H_{2}:=\left\{(x,y,\varphi,\psi)\in H|y\leq x, y\leq 0\right\}\]
\[H_{3}:=\left\{(x,y,\varphi,\psi)\in H|0\leq x, 0\leq y\right\}.\]
\end{defn}
This subdivision of the complex line can be found in $\cite{Vaf}$.
The inequalities defining the subdivision of $H$ determine, in particular,  a subdivision in three parts of the amoeba $\pi(H)$. This subdivision is represented in the following picture:

\begin{center}
\includegraphics[scale=0.6]{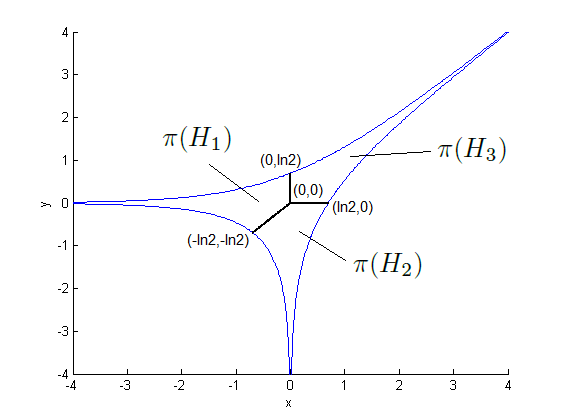}
\begin{center}
\textbf{Figure 5}
\end{center}
\end{center}
Analogously, taking the image under the argument map of $H_{1}$, $H_{2}$ and $H_{3}$ we get an induced subdivision of $Arg(H)$ in three parts.
\begin{defn}\label{defn:5}
We set: 
\begin{alignat*}{9}
A_{1}&:=\left\{(\varphi,\psi)\in Arg(H_{1})|0\leq \varphi \leq\pi\right\},\quad &A_{2}&:=\left\{(\varphi,\psi)\in Arg(H_{1})|\pi\leq \varphi \leq2\pi\right\},\\
B_{1}&:=\left\{(\varphi,\psi)\in Arg(H_{2})|0\leq \varphi\leq \pi\right\},\quad &B_{2}&:=\left\{(\varphi,\psi)\in Arg(H_{2})|\pi\leq \varphi\leq 2\pi\right\}, \\
C_{1}&:=\left\{(\varphi,\psi)\in Arg(H_{3})|0\leq \varphi \leq\pi\right\},\quad &C_{2}&:=\left\{(\varphi,\psi)\in Arg(H_{3})|\pi\leq \varphi \leq2\pi\right\}. 
\end{alignat*}
\end{defn}

\begin{prop}
\begin{alignat*}{9}
A_{1}&=\left\{(\varphi,\psi)\in Arg(H)|\psi\leq\frac{\varphi}{2}+\pi,\psi\leq-\varphi+2\pi\right\},\\
A_{2}&=\left\{(\varphi,\psi)\in Arg(H)|\psi\geq\frac{\varphi}{2},\psi\geq-\varphi+2\pi\right\}.
\end{alignat*}
\end{prop}
\proof
From Proposition $\ref{prop:2}$ we immediately see that $\left\{(0,\pi),(\pi,\pi)\right\}$ are both elements of $A_{1}$ and $A_{2}$.
Using Lemma $\ref{lem:2}$ we get that $H_{1}\cap H_{\circ}$ equals
\[\left\{(x,y,\varphi,\psi)\in H|x\leq y,x\leq0, e^{x}=-\frac{\sin\psi}{\sin(\psi-\varphi)},e^{y}=\frac{\sin\varphi}{\sin(\psi-\varphi)}\right\}\]
that is, $H_{1}\cap H_{\circ}$ equals:
\[\left\{(x,y,\varphi,\psi)\in H|-\frac{\sin\psi}{\sin(\psi-\varphi)}\leq\frac{\sin\varphi}{\sin(\psi-\varphi)},-\frac{\sin\psi}{\sin(\psi-\varphi)}\leq 1\right\}.\]
Therefore, $Arg(H_{1})\cap Arg(H_{\circ})$ is equal to:
\begin{equation}\label{eq:r1}
\left\{(\varphi,\psi)\in Arg(H)|-\frac{\sin\psi}{\sin(\psi-\varphi)}\leq\frac{\sin\varphi}{\sin(\psi-\varphi)},-\frac{\sin\psi}{\sin(\psi-\varphi)}\leq 1\right\}.
\end{equation}
By Definition $\ref{defn:5}$, $A_{1}$ is given by the subset of $Arg(H_{1})$ characterized by $0\leq\varphi\leq \pi$. Thus, let us assume $\varphi\in\left(0,\pi\right)$, then $\psi\in(\pi,\varphi+\pi)$ and hence $\sin(\psi-\varphi)>0$. Therefore, from $\eqref{eq:r1}$ we get that the subset of $Arg(H_{1})\cap Arg(H_{\circ})$ such that $\varphi\in(0,\pi)$ is equal to:
\begin{equation}\label{eq:r2}
\left\{(\varphi,\psi)\in Arg(H)|0<\varphi<\pi,-\sin\psi\leq\sin\varphi,-\sin\psi\leq\sin(\psi-\varphi)\right\},
\end{equation}
which, in turn, equals the subset of $(S^{1})^{2}$ satisfying:
\begin{equation}\label{eq:s1}
\begin{cases}
0<\varphi<\pi\\
\pi<\psi<\varphi+\pi\\
-\sin\psi\leq\sin\varphi\\
-\sin\psi\leq\sin(\psi-\varphi)
\end{cases}.
\end{equation}
Now, if $\varphi\left(0,\frac{\pi}{2}\right]$ the solution of: 
\begin{equation}\label{eq:s2}
\begin{cases}
0<\varphi<\pi\\
\pi<\psi<\varphi+\pi\\
-\sin\psi\leq\sin\varphi\\
\end{cases}
\end{equation}
is given by all $(\varphi,\psi)\in Arg(H)$ such that $\varphi\in\left(0,\frac{\pi}{2}\right]$. If $\varphi\in\left[\frac{\pi}{2},\pi\right)$, the solution of $\eqref{eq:s2}$ is given by all $(\varphi,\psi)\in Arg(H)$ such that $\psi\leq-\varphi+2\pi$. 
On the other hand, the solution of:
\begin{equation}\label{eq:o}
\begin{cases}
0<\varphi<\pi\\
\pi<\psi<\varphi+\pi\\
-\sin\psi\leq\sin(\psi-\varphi)
\end{cases}
\end{equation}
is given by all $(\varphi,\psi)\in Arg(H)$ such that $\psi\leq\frac{\varphi}{2}+\pi$. Indeed, let us fix $\psi\in(\pi,2\pi)$. Using the first two inequalities in $\eqref{eq:o}$ we get $\psi-\pi<\psi-\varphi<\pi$. Now, let $\psi\in\left(\frac{3\pi}{2},2\pi\right]$. Then, the third inequality in $\eqref{eq:o}$ is satisfied if and only if $2\pi-\psi\leq\psi-\varphi\leq\psi-\pi$. But, the system of inequalities:
\[
\begin{cases}
0<\varphi<\pi\\
\psi-\pi<\psi-\varphi<\pi\\2\pi-\psi\leq\psi-\varphi\leq\psi-\pi
\end{cases}
\]
has no solution. If $\psi\in\left(\pi,\frac{3\pi}{2}\right]$, then the third inequality in $\eqref{eq:o}$ is satisfied if and only if $\psi-\pi\leq\psi-\varphi\leq2\pi-\psi$. In this case, the solution of the system of inequalities:
\[
\begin{cases}
0<\varphi<\pi\\
\psi-\pi<\psi-\varphi<\pi\\
\psi-\pi\leq\psi-\varphi\leq2\pi-\psi
\end{cases}
\] 
is given by all $(\varphi,\psi)\in Arg(H)$ such that $\psi\leq\frac{\varphi}{2}+\pi$ , as claimed. Thus, the statement for $A_{1}$ is proved.
By Definition $\ref{defn:5}$, $A_{2}$ is given by the subset of $Arg(H_{1})$ for $\pi\leq\varphi\leq 2\pi$. Thus, let us assume $\varphi\in\left(\pi,2\pi\right)$, then $\psi\in(0,\pi)$ and hence $\sin(\psi-\varphi)<0$. Therefore, from $\eqref{eq:r1}$ we get that the subset of $Arg(H_{1})\cap Arg(H_{\circ})$ such that $\varphi\in(\pi,2\pi)$ is equal to:
\begin{equation}\label{eq:r3}
\left\{(\varphi,\psi)\in Arg(H)|\pi<\varphi<2\pi,-\sin\psi\geq\sin\varphi,-\sin\psi\geq\sin(\psi-\varphi)\right\}.
\end{equation}
In the same way as for $A_{1}$, one can check that the inequalities in $\eqref{eq:r3}$ induce the inequalities in the statement for $A_{2}$.

\endproof

In the same way, one can prove the following
\begin{prop}
\begin{alignat*}{9}
B_{1}&=\left\{(\varphi,\psi)\in Arg(H)|\psi\leq2\varphi,\psi\geq-\varphi+2\pi\right\},\\
B_{2}&=\left\{(\varphi,\psi)\in Arg(H)|\psi\geq2(\varphi-\pi),\psi\leq-\varphi+2\pi\right\}.
\end{alignat*}
\end{prop}

\begin{prop}
\begin{alignat*}{9}
C_{1}&=\left\{(\varphi,\psi)\in Arg(H)|\psi\geq\frac{\varphi}{2}+\pi,\psi\geq2\varphi\right\},\\
C_{2}&=\left\{(\varphi,\psi)\in Arg(H)|\psi\leq2(\varphi-\pi),\psi\leq\frac{\varphi}{2}\right\}.
\end{alignat*}
\end{prop}

The subdivision of $Arg(H)$ is represented in the following picture:
\begin{center}
\includegraphics[scale=0.6]{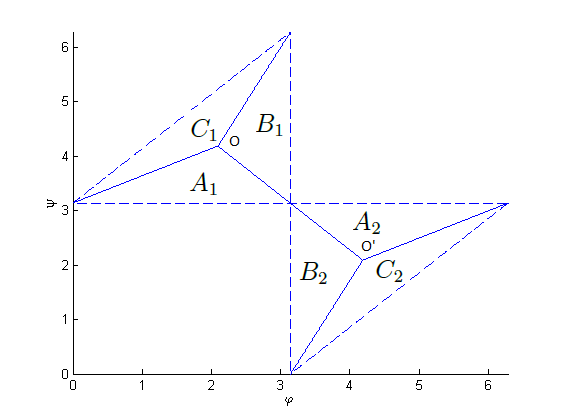}
\begin{center}
\textbf{Figure 6}
\end{center}
\end{center}
The two points $O:=(\frac{2\pi}{3},\frac{4\pi}{3})$ and $O':=(\frac{4\pi}{3},\frac{2\pi}{3})$ are given by the intersection of $\psi=-\varphi+2\pi$ respectively with $\psi=\frac{\varphi}{2}$ and $\psi=2\varphi$. 
These two points are the arguments of the fibre of $\pi|_{H}$ over $(0,0)$. Indeed, using Lemma $\ref{lem:1}$ we get that the fibre of $\pi|_{H}$ over $(0,0)$ is given by $\left\{(0,0,\varphi,\psi),(0,0,2\pi-\varphi,2\pi-\psi)\right\}$ or $\left\{(0,0,2\pi-\varphi,\psi),(0,0,\varphi,2\pi-\psi)\right\}$, where $(0,0,\varphi,\psi)$ is a solution of the following system of equations: 
\[
\begin{cases}
e^{2x}=1+2e^{y}\cos\psi+e^{2y}\\e^{2y}=1+2e^{x}\cos\varphi+e^{2x}
\end{cases},
\]
that is:
\[
\begin{cases}
2\cos\psi+1=0\\2\cos\varphi+1=0
\end{cases}.
\]
Take $(\varphi,\psi)=(\frac{2\pi}{3},\frac{4\pi}{3})$, then we get that the two above sets coincide and $(2\pi-\varphi,2\pi-\psi)=(\frac{4\pi}{3},\frac{2\pi}{3})$.
We notice that $O$ and $O'$ are the barycentres respectively of the triangles $\overline{T_{1}}$ and $\overline{T_{2}}$, closure in $(S^{1})^{2}$ of the subsets $T_{1}$ and $T_{2}$.

In a similar way as done for $H$, we subdivide $H_{trop}$ in three parts, which we call $H_{1trop}$, $H_{2trop}$, $H_{3trop}$. To do that, we extend the above subdivision of $Arg(H)$ to $\overline{Arg(H)}$:
\begin{center}
\includegraphics[scale=0.6]{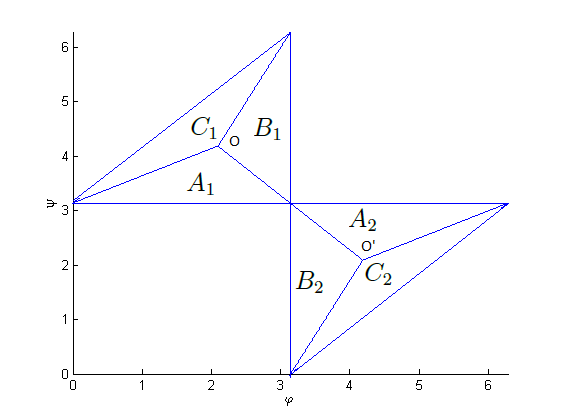}
\begin{center}
\textbf{Figure 7}
\end{center}
\end{center}
\begin{defn}\label{defn:44}
We set:
\[H_{1trop}:=\left\{(x,0,\varphi,\psi)\in H_{trop}|(\varphi,\psi)\in \overline{Arg(H_{1})}\right\},\]
\[H_{2trop}:=\left\{(0,y,\varphi,\psi)\in H_{trop}|(\varphi,\psi)\in \overline{Arg(H_{2})}\right\},\]
\[H_{3trop}:=\left\{(x,x,\varphi,\psi)\in H_{trop}|(\varphi,\psi)\in \overline{Arg(H_{3})}\right\}.\]
\end{defn}

\begin{lem}\label{lem:4}
Let us consider the automorphism of $\mathbb{R}^{2}\times (S^{1})^{2}$ of order $3$:
\begin{equation}
\begin{split}
\lambda:\mathbb{R}^{2}\times (S^{1})^{2}&\longrightarrow \mathbb{R}^{2}\times (S^{1})^{2}\\
(x,y,\varphi,\psi)&\longmapsto (-y,x-y,-\psi+2\pi,\varphi-\psi+2\pi)\\
\end{split}.
\end{equation}
 We have $\lambda(H_{1})=H_{2}$, $\lambda(H_{2})=H_{3}$, $\lambda(H_{3})=H_{1}$ and $\lambda(H_{1trop})=H_{2trop}$, $\lambda(H_{2trop})=H_{3trop}$, $\lambda(H_{3trop})=H_{1trop}$.
\end{lem}
\proof
We only show $\lambda(H_{1})=H_{2}$ and $\lambda(H_{1trop})=H_{2trop}$. Let $P:=(x,y,\varphi,\psi)\in H$. Then $P\in H_{1}$ if and only if it satisfies:
\begin{equation}\label{eq:l1}
\begin{cases}
e^{x+i\varphi}+e^{y+i\psi}+1=0 \\ x\leq 0\\y\geq x
\end{cases}.
\end{equation}
Similarly, $P\in H_{2}$ if and only if it satisfies:
\begin{equation}\label{eq:l3}
\begin{cases}
e^{x+i\varphi}+e^{y+i\psi}+1=0 \\ y\leq x\\y\leq 0
\end{cases}.
\end{equation}
Assume $P\in H_{1}$. Then $\lambda(P)=(-y,x-y,-\psi+2\pi,\varphi-\psi+2\pi)$ satisfies:  
\begin{equation}\label{eq:l2}
\begin{cases}
e^{-y+i(-\psi)}+e^{x-y+i(\varphi-\psi)}+1=0 \\ x-y\leq -y \\ x-y\leq 0
\end{cases}.
\end{equation}
Indeed, $e^{-y+i(-\psi+2\pi)}+e^{x-y+i(\varphi-\psi+2\pi)}+1=( 1+e^{x+i\varphi}+e^{y+i\psi})(e^{y+i\psi})^{-1}$. Using the first equation in $\eqref{eq:l1}$, we get that the first equation in $\eqref{eq:l2}$ holds. Moreover, the first and second inequality in $\eqref{eq:l2}$ are direct consequences respectively of the first and second inequality in $\eqref{eq:l1}$. Thus, changing coordinates, $\lambda(P)=:(x^{'},y^{'},\varphi^{'},\psi^{'})$, we see that $\lambda(P)\in H_{2}$.
Conversely, let $P':=(x,y,\varphi,\psi)$ satisfy $\eqref{eq:l3}$. Then, in a similar way as before one shows that $\lambda^{-1}(P')=(y-x,-x,\psi-\varphi+2\pi,-\psi+2\pi)$ satisfies $\eqref{eq:l1}$. To show that $\lambda(H_{1trop})=H_{2trop}$ we notice that from $\lambda(H_{1})=H_{2}$, it immediately follows that $\lambda$ maps $\pi(H_{1trop})$ onto $\pi(H_{2trop})$ and $Arg(H_{1trop})=\overline{Arg(H_{1})}$ onto $Arg(H_{2trop})=\overline{Arg(H_{2})}$. 
\endproof 

\begin{rem}
The automorphism $\lambda$ of $\mathbb{R}^{2}\times (S^{1})^{2}$ descends from the following automorphism of $(\mathbb{C}^{*})^{2}$ of order $3$:
\begin{equation}
\begin{split}
\lambda^{'}:(\mathbb{C}^{*})^{2}&\longrightarrow(\mathbb{C}^{*})^{2}\\
(z,w)&\longmapsto (w^{-1},zw^{-1})\\
\end{split}
\end{equation}
via the identification 
\begin{equation}
\begin{split}
h:(\mathbb{C}^{*})^{2}&\longrightarrow \mathbb{R}^{2}\times (S^{1})^{2}\\
(z_1,z_2)&\longmapsto(\ln|z_1|,\ln|z_2|,\arg z_1,\arg z_2)\\
\end{split}.
\end{equation}
\end{rem}

We further subdivide $H_{1}$ and $H_{1trop}$ in two parts.

\begin{defn}
We set:
\[ H_{1L}:=\left\{(x,y,\varphi,\psi)\in H_{1}|y\geq 2x+\ln2\right\}\]
and
\[ H_{1T}:=\left\{(x,y,\varphi,\psi)\in H_{1}|y\leq 2x+\ln2\right\}.\]
\end{defn}
Taking the image under $\pi$ of $H_{1L}$ and $H_{1T}$ we get an induced subdivision of $\pi(H_{1})$. 

\begin{defn}
We set:
\[\Gamma_{1}:=\left\{(x,y,\varphi,\psi)\in H_{1}|y=2x+\ln2\right\}.\]
\end{defn}

The subdivision of $\pi(H_{1})$ is represented in the following picture:
\begin{center}
\includegraphics[scale=0.6]{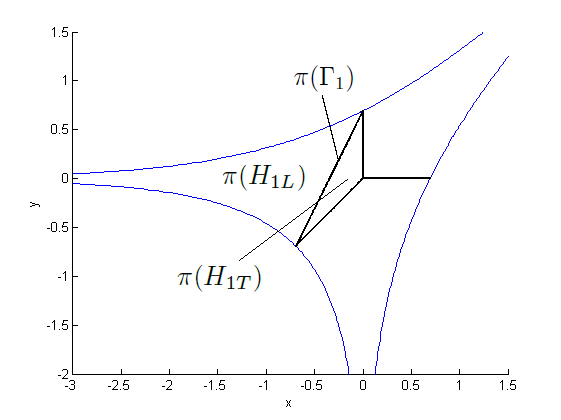}
\begin{center}
\textbf{Figure 8}
\end{center}
\end{center}
The picture justifies the names given to the two parts of the subdivision. Namely, $'L'$ stands for $'Leg'$ and $'T'$ for $'Triangle'$.

Taking the image under $Arg$ of $H_{1L}$ and $H_{1T}$ we get an induced subdivision of $Arg(H_{1})$. Using Proposition $\ref{prop:2}$ and Lemma $\ref{lem:2}$ one can immediately prove:
\begin{prop}\label{prop:31}
$Arg(\Gamma_{1})$ equals:
\[\left\{(\varphi,\psi)\in Arg(H)| \sin^{2}\psi=\frac{1}{2}\sin\varphi\sin(\psi-\varphi),-\frac{\sin\psi}{\sin(\psi-\varphi)}> 0\right\}\cup\left\{(0,\pi),(\pi,\pi)\right\}.\]
\end{prop}
\proof
Using Proposition $\ref{prop:2}$ we immediately get that $\left\{(0,\pi),(\pi,\pi)\right\}\in Arg(\Gamma_{1})$. 
By Lemma $\ref{lem:2}$, we get that $\Gamma_{1}\cap H_{\circ}$ equals
\[\left\{(x,y,\varphi,\psi)\in H|y=2x+\ln2, e^{x}=-\frac{\sin\psi}{\sin(\psi-\varphi)},e^{y}=\frac{\sin\varphi}{\sin(\psi-\varphi)}\right\}\]
that is, $\Gamma_{1}\cap H_{\circ}$ equals:
\[\left\{(x,y,\varphi,\psi)\in H|(-\frac{\sin\psi}{\sin(\psi-\varphi)})^{2}=\frac{1}{2}\frac{\sin\varphi}{\sin(\psi-\varphi)},-\frac{\sin\psi}{\sin(\psi-\varphi)}> 0\right\}.\]
Therefore, $Arg(\Gamma_{1})\cap Arg(H_{\circ})$ is equal to:
\begin{equation}\label{eq:g}
\left\{(\varphi,\psi)\in Arg(H)|\sin^{2}\psi=\frac{1}{2}\sin\varphi\sin(\psi-\varphi),-\frac{\sin\psi}{\sin(\psi-\varphi)}> 0\right\}.
\end{equation}
\endproof
The subdivision of $Arg(H_{1})$ is represented in the following picture:
\begin{center}
\includegraphics[scale=0.7]{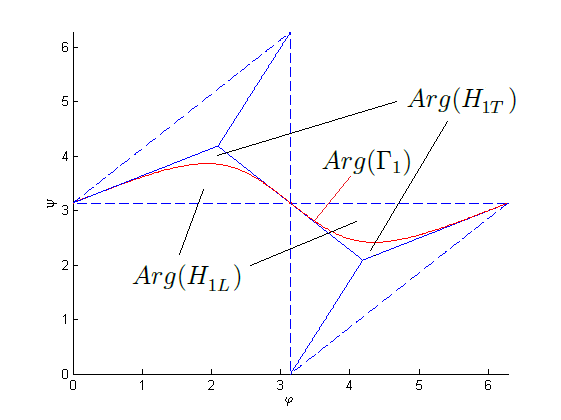}
\begin{center}
\textbf{Figure 9}
\end{center}
\end{center}

\begin{lem}\label{lem:9}
The map:
\[Arg|_{H_{1T}}:H_{1T}\longrightarrow Arg(H_{1T})\]
is a homeomorphism.
\end{lem}
\proof
By Proposition $\ref{prop:3}$ we get that $Arg|_{H_{1T}}$ is bijective on $H_{1T}\cap H_{\circ}$. By Proposition $\ref{prop:2}$, we see that if $(\varphi,\psi)\in\left\{(0,\pi),(\pi,\pi)\right\}$, then there exists only one point in $H_{1T}$ which is mapped to $(\varphi,\psi)$ under $Arg$. Moreover, it is continuous being the restriction of a continuous map. We easily see that $H_{1T}$ is compact and since $Arg(H_{1T})$ is Hausdorff, we get that $Arg|_{H_{1T}}$ is a homeomorphism.
\endproof

\begin{defn}\label{defn:45}
We set:
\[H_{1tropL}:=\left\{(x,0,\varphi,\pi)\in H_{1trop}|x\leq 0\right\},  H_{1tropT}:=\left\{(0,0,\varphi,\psi)\in H_{1trop}\right\}.\]
\end{defn}

\begin{lem}\label{lem:5}
Let $c\geq\ln2$. The straight line $r\subset\mathbb{R}^{2}$ of equation $y=2x+c$ intersects the boundary of the amoeba of $H$ in two points. Indeed, it has no intersection points with the curve $e^{x}-e^{y}=1$. The intersection point between $r$ and the curve $e^{y}-e^{x}=1$ is unique and it is given by $x=\ln(\frac{1+\sqrt{1+4e^{c}}}{2e^{c}})$. Moreover, $r$ intersects the curve $e^{y}+e^{x}=1$ in a unique point given by $x=\ln(\frac{-1+\sqrt{1+4e^{c}}}{2e^{c}})$.
\end{lem}
\proof
The system of equations:
\[
\begin{cases}
e^{x}-e^{y}=1\\
y=2x+c
\end{cases}
\]
is equivalent to:
\begin{equation}\label{eq:l1a}
\begin{cases}
e^{c}e^{2x}-e^{x}+1=0\\ y=2x+c
\end{cases}.
\end{equation}
Solving the first equation in $\eqref{eq:l1a}$ in the variable $e^{x}$, we see that it has no solution in $\mathbb{R}$.
The system of equations:
\[
\begin{cases}
e^{y}-e^{x}=1\\
y=2x+c
\end{cases}
\]
is equivalent to:
\begin{equation}\label{eq:l1b}
\begin{cases}
e^{c}e^{2x}-e^{x}-1=0\\ y=2x+c
\end{cases}.
\end{equation}
Solving the first equation in $\eqref{eq:l1b}$ in the variable $e^{x}$ we see that it gives two solutions, namely $e^{x}=\frac{1-\sqrt{1+4e^{c}}}{2e^{c}}$ and $e^{x}=\frac{1+\sqrt{1+4e^{c}}}{2e^{c}}$. Since $c\geq\ln2$, we have $e^{c}\geq 2$ and hence $\sqrt{1+4e^{c}}\geq 3$. Thus, $\frac{1-\sqrt{1+4e^{c}}}{2e^{c}}<0$.  Therefore, the solution of $\eqref{eq:l1b}$ is unique and it is given by $x=\ln(\frac{1+\sqrt{1+4e^{c}}}{2e^{c}})$.
In the same way, we get that the system of equations:
\[
\begin{cases}
e^{x}+e^{y}=1\\
y=2x+c
\end{cases}
\]
has a unique solution and it is given by $x=\ln(\frac{-1+\sqrt{1+4e^{c}}}{2e^{c}})$.
\endproof

\begin{lem}\label{lem:6}
Let $a,b\in\mathbb{R}$, $a<b$ and $f:\left[a,b\right]\longrightarrow\mathbb{R}$ be a differentiable strictly concave function. Let $P\in \mathbb{R}^{2}$ and $c\in\left[a,b\right]$. Let $r:=\left\{(x,y)\in\mathbb{R}^{2}|y=mx+q\right\}$ the straight line passing through $P$ and through $(c,f(c))$. If $m\geq f^{'}(a)$ or $m\leq f^{'}(b)$ , then $(c,f(c))\in\mathbb{R}^{2}$ is the unique point of intersection between $r$ and the graph of $f$.
\end{lem}
\proof
Assume $m\geq f^{'}(a)$.
By contradiction, let $(d,f(d))\in\mathbb{R}^{2}$ another point of intersection, with $d\in \left[a,b\right]$ and assume $d>c$. Then, by the Mean Value Theorem there exists $c^{'}\in(c,d)$ such that $f^{'}(c^{'})=m$. Therefore, we get $f^{'}(c^{'})\geq f^{'}(a)$. On the other hand, since the derivative of a strictly concave function is monotonically strictly decreasing, being $a<c^{'}$, we get $f^{'}(a)> f^{'}(c^{'})$.
A similar contradiction we get assuming $m\leq f^{'}(b)$. 
\endproof

Let $c\geq\ln2$ and consider the subset of $H_{1L}$:
\[W_{c}:=\left\{(x,y,\varphi,\psi)\in H_{1L}|y=2x+c\right\}.\] We have:
\[H_{1L}=\bigsqcup_{c\geq\ln2}W_{c}.\] The next proposition shows that if $(\varphi_0,\psi_0)\in Arg(W_c)\subset Arg(H_{1L})$, then the straight line passing through $O:=(\frac{2\pi}{3},\frac{4\pi}{3})\in Arg(H)$ and $(\varphi_0,\psi_0)$, for $0\leq\varphi_{0}\leq \pi$, and  through $O^{'}:=(\frac{4\pi}{3},\frac{2\pi}{3})\in Arg(H)$ and $(\varphi_0,\psi_0)$, for $\pi\leq\varphi_{0}\leq 2\pi$, does not have any other point of intersection with $Arg(W_c)$ apart from $(\varphi_0,\psi_0)$ itself. This statement rigorously proves what is visually evident from the pictorial representation of $Arg(W_c)$ (see Figure 10 and Figure 11 below).  
\begin{prop}\label{prop:4}
Let $c\geq \ln 2$ and $W_{c}:=\left\{(x,y,\varphi,\psi)\in H|y=2x+c\right\}$. Let us consider a point  $(x_{0},y_{0},\varphi_{0},\psi_{0})\in~ W_{c}$, with $0\leq\varphi_{0}\leq \pi$. Let $s$ be the straight line in $(S^{1})^{2}$ passing through $(\varphi_{0},\psi_{0})$ and $O:=(\frac{2\pi}{3},\frac{4\pi}{3})$. Then, for $0\leq\varphi\leq \pi$, $(\varphi_{0},\psi_{0})$ is the unique point of intersection between $s$ and $Arg(W_{c})$. Similarly, let $(x_{0},y_{0},\varphi_{0},\psi_{0})\in W_{c}$, with $\pi\leq\varphi_{0}\leq 2\pi$. Let $s^{'}$ be the straight line in $(S^{1})^{2}$ passing through $(\varphi_{0},\psi_{0})$ and $O^{'}:=(\frac{4\pi}{3},\frac{2\pi}{3})$. Then, for $\pi\leq\varphi\leq 2\pi$, $(\varphi_{0},\psi_{0})$ is the unique point of intersection between $s^{'}$ and $Arg(W_{c})$.
\end{prop}
\proof
By Proposition $\ref{prop:2}$ and Lemma $\ref{lem:2}$, we have that $(\varphi,\psi)\in Arg(W_{c})$ if and only if either $(\varphi,\psi)\in\left\{(0,\pi),(\pi,\pi)\right\}$ or it satisfies:
\begin{equation}\label{eq:s4}
\begin{cases}
e^{x}=-\frac{\sin\psi}{\sin(\psi-\varphi)}
\\ e^{y}=+\frac{\sin\varphi}{\sin(\psi-\varphi)}\\
y=2x+c
\end{cases}
\end{equation}
which is equivalent to:
\begin{equation} \label{eq:s5}
\begin{cases}
 \sin^{2}(\psi)=k \sin(\varphi)\sin(\psi-\varphi)\\-\frac{\sin\psi}{\sin(\psi-\varphi)}>0\\
+\frac{\sin\varphi}{\sin(\psi-\varphi)}>0
\end{cases}
\end{equation} 
where $k:=\frac{1}{e^{c}}$. 
We show that $Arg(W_{c})$, seen as a function, $\psi=f(\varphi)$, is a differentiable strictly concave function on the interval $0\leq\varphi\leq\pi$. By implicit differentiation, we get that its first derivative is:
\begin{equation}\label{eq:s6}
\frac{d\psi}{d\varphi}=\frac{2k\sin(\psi-2\varphi)}{2\sin(2\psi)-k\sin\psi+k\sin(\psi-2\varphi)}.
\end{equation}
We have that $\frac{d\psi}{d\varphi}$ is continuous on $0<\varphi<\pi$ and it can be naturally extended to obtain continuity on the closed interval.
Moreover, we have that $\frac{d\psi}{d\varphi}$ is strictly decreasing. To show that, let $\pi_{+}$ be as in Proposition $\ref{prop:33}$, we notice that the curve $\eqref{eq:s5}$  corresponds in $\pi(H)$ to the path given by the line segment $V:=\left\{y=2x+c\right\}\cap\pi(H)$ via the diffeomorphism $Arg\circ\pi^{-1}_{+}|_{V}$, which maps $(0,\pi),(\pi,\pi)\in Arg(W_{c})$ to the points in $V$ respectively with coordinate $x_{0}$ and $x_{1}$, where:
\[e^{x_{0}}=\frac{k}{2}(-1+\sqrt{1+\frac{4}{k}}) \text{ and } e^{x_{1}}=\frac{k}{2}(1+\sqrt{1+\frac{4}{k}}),\]
as stated in Lemma $\ref{lem:5}$. Therefore, we want to rewrite $\eqref{eq:s6}$ in the $(x,y)$ coordinates. We have that $\eqref{eq:s6}$ equals:

\[
\frac{2k(\cos(2\varphi)-2\cos\varphi\cos\psi\frac{\sin\varphi}{\sin\psi})}{4\cos\psi-k+k(\cos(2\varphi)-2\cos\varphi\cos\psi\frac{\sin\varphi}{\sin\psi})}.
\]

Using $\eqref{eq:s4}$, the identity $\cos(2\varphi)=-1+2\cos^{2}\varphi$ and the expressions for $\cos\varphi,\cos\psi$ given in $\eqref{eq:1}$, we can rewrite the last expression as:
\[\frac{2k(-e^{2y}+1)}{2(e^{4x-y}-e^{2x+y}-e^{2x-y})-k(e^{2x}+e^{2y}-1)},\]
which, using the relations $y=2x+c$ and $k=\frac{1}{e^{c}}$, can be rewritten as:
\begin{equation}\label{eq:s10}
\frac{2e^{4x}-2k^{2}}{3e^{4x}-k^{2}e^{2x}+k^{2}}.
\end{equation}
Now,
\begin{equation}\label{eq:s11}
\frac{d}{dx}(\frac{2e^{4x}-2k^{2}}{3e^{4x}-k^{2}e^{2x}+k^{2}})=4k^{2}e^{2x}\frac{-e^{4x}+8e^{2x}-k^{2}}{(3e^{4x}-k^{2}e^{2x}+k^{2})^{2}}.
\end{equation}
 We have that $\eqref{eq:s11}$ is positive if and only if $4-\sqrt{16-k^{2}}<e^{2x}<4+\sqrt{16-k^{2}}$. For $x\in\left[x_{0},x_{1}\right]$, we have that $\eqref{eq:s11}$ is always positive. Moreover, by $\eqref{eq:1}$ we have:
\[\cos\varphi=\frac{\frac{e^{4x}}{k^{2}}-e^{2x}-1}{2e^{x}},\]
hence for $0\leq\varphi\leq\pi$, if $\varphi$ strictly increases, then $x$ strictly decreases. Thus, $\frac{d\psi}{d\varphi}$ is strictly decreasing. Now, we have that the straight line passing through $O$ and $(0,\pi)$ has angular coefficient equal to $\frac{1}{2}$, while the straight line passing through $O$ and $(\pi,\pi)$ has angular coefficient equal to $-1$. Therefore, let $(\varphi_{0},\psi_{0})\in Arg(W_{c})$, with $0\leq\varphi_{0}\leq\pi$, and $s$ be as in the statement. One easily sees that the angular coefficient $m$ of $s$ satisfies $m\geq \frac{1}{2}$ or $m\leq -1$. On the other hand, we notice that if $c\geq \ln 2$, then $k\leq\frac{1}{2}$, which, in turn, implies  $\frac{d\psi}{d\varphi}((0,\pi))\leq\frac{1}{2}$. Similarly, $k\leq\frac{1}{2}$ implies $\frac{d\psi}{d\varphi}((\pi,\pi))\geq -1$.
Therefore, since by assumption $c\geq \ln2$, we have $m\geq\frac{d\psi}{d\varphi}((0,\pi))$ or $m\leq\frac{d\psi}{d\varphi}((\pi,\pi))$. Thus, by Lemma $\ref{lem:6}$, we have that it does not exist another point of intersection $(\varphi_{1},\psi_{1})$ between $s$ and $Arg(W_{c})$, with $0\leq\varphi_{1}\leq\pi$.
The second part of the statement follows from the fact that $H$ is symmetric with respect to the transformation $(x,y,\varphi,\psi)\longmapsto (x,y,2\pi-\varphi,2\pi-\psi)$.
\endproof

\begin{center}
\includegraphics[scale=0.6]{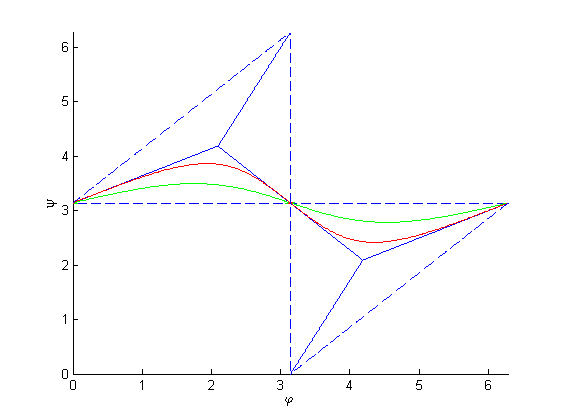}

\textbf{Figure 10}: In red $Arg(W_{c})$ for $k=\frac{1}{2}$, in green for $k=\frac{1}{8}$.
\end{center}

\section{The isotopy map}
\label{sec: isotopy}

\begin{defn}
Let $A_{1},A_{2}\subset Arg(H_{1})$ as in Definition $\ref{defn:5}$ and let $(\varphi,\psi)\in A_{1}$. Let $Arg(\Gamma_{1})$ be as in Proposition $\ref{prop:31}$, let $O:=(\frac{2\pi}{3},\frac{4\pi}{3})\in A_1$ and $r\subset (S^{1})^{2}$ the straight line passing through $O$ and through $(\varphi,\psi)$. Let $Q'\in\overline{Arg(H_{1})}$ the point of intersection between $r$ and the straight line in $(S^{1})^{2}$ of equation $\psi=\pi$.\\
Let $d:(S^{1})^{2}\times (S^{1})^{2}\rightarrow \mathbb{R}_{\geq 0}$ be the flat metric on $(S^{1})^{2}$. We set:
\[b:=d(O,Q').\]
If $(\varphi,\psi)\in A_{1}\cap Arg(H_{1T})$, let $Q$ be the intersection point between $r$ and $Arg(\Gamma_{1})$, we set:
\[a:=d(O,Q).\]
If $(\varphi,\psi)\in A_{1}\cap Arg(H_{1L})$, then we set:
\[a':=d(O,(\varphi,\psi)).\]
See Figure 11 below:

\begin{center}
\includegraphics[scale=0.8]{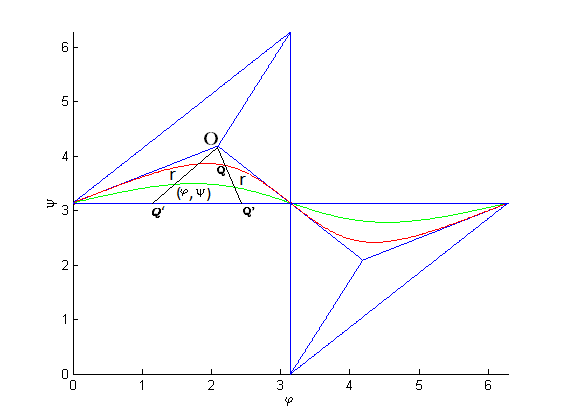}
\begin{center}
\textbf{Figure 11}
\end{center}
\end{center}

Now, for any $t\in\left[0,1\right]$ we define a map:
\[
(\varphi_{t},\psi_{t}): Arg(H_{1})\longrightarrow \overline{Arg(H_{1})},
\]
by:
\begin{equation}\label{eq:t1}
(\varphi_{t},\psi_{t})|_{A_{1}}:(\varphi,\psi)\longmapsto 
\begin{cases}
O+(\frac{b}{a})^{t}((\varphi,\psi)-O),&\text{if } (\varphi,\psi)\in A_{1}\cap Arg(H_{1T})\setminus\left\{O\right\};\\&\\
O+(\frac{b}{a'})^{t}((\varphi,\psi)-O),&\text{if } (\varphi,\psi)\in A_{1}\cap Arg(H_{1L})\setminus\left\{O\right\};\\&\\
O,&\text{if } (\varphi,\psi)= O.
\end{cases}
\end{equation}
The definition of $(\varphi_{t},\psi_{t})|_{A_{2}}$ is obtained by replacing $O$ with $O':=(\frac{4\pi}{3},\frac{2\pi}{3})$ and $A_{1}$ with $A_{2}$ in the above definitions of $a, a', b$ and $(\varphi_{t},\psi_{t})|_{A_{1}}$. 
\end{defn}

\begin{claim}\label{cl:7}
The map $(\varphi_{t},\psi_{t})$ is well defined.
\end{claim}
\proof
If $(\varphi,\psi)\in A_{1}\cap Arg(H_{1T})$, then $a$ is well defined as the intersection point $Q$ between $r$ and $Arg(\Gamma_{1})$ is unique by Proposition $\ref{prop:4}$.
We immediately see that on the intersection between $Arg(H_{1T})$ and $Arg(H_{1L})$, that is on the curve $Arg(\Gamma_{1})$, the definitions of $a$ and $a'$ coincide. Moreover, on $A_{1}\cap A_{2}$, $(\varphi_{t},\psi_{t})|_{A_{1}}=(\varphi_{t},\psi_{t})|_{A_{2}}$.
\endproof

\begin{lem}\label{lem:8}
The map $(\varphi_{t},\psi_{t})$ is continuous. Moreover, it is injective on $Arg(H_{1T})$ and $(\varphi_{1},\psi_{1})(Arg(H_{1T}))=\overline{Arg(H_{1})}$.
\end{lem}
\proof
The map is clearly continuous. From the definition of $(\varphi_{t},\psi_{t})$ one easily sees that $O$ and $O^{'}$ are the unique points respectively sent to $O$ and $O^{'}$. Now, let $(\varphi,\psi)\neq O,O'$ and $(\varphi',\psi')\neq O,O'$ in $Arg(H_{1T})$. By construction, we have that if $(\varphi_{t},\psi_{t})(\varphi,\psi)=(\varphi_{t},\psi_{t})(\varphi',\psi')$ , then $(\varphi,\psi)$ and $(\varphi',\psi')$ lie on the same straight line passing through $O$ or through $O'$. Since this straight line uniquely determines $a$ and $b$, we get $(\varphi,\psi)=(\varphi',\psi')$.  Finally, the last part of the statement follows immediately from the definition of $(\varphi_{1},\psi_{1})$.
\endproof

\begin{defn}\label{defn:6}
 Let $P=(x,y,\varphi,\psi)\in H_{1}$. For each $t\in\left[0,1\right]$, we define maps:
\begin{equation}
\begin{split}
\Phi^{t}_{1T}:H_{1T}&\longrightarrow \mathbb{R}^{2}\times (S^{1})^{2}\\
(x,y,\varphi,\psi)&\longmapsto (x(1-t),y(1-t),\varphi_{t},\psi_{t})\\
\end{split}
\end{equation}
and
\begin{equation}
\begin{split}
\Phi^{t}_{1L}:H_{1L}&\longrightarrow \mathbb{R}^{2}\times (S^{1})^{2}\\
(x,y,\varphi,\psi)&\longmapsto (x-t\frac{y-\ln2}{2},y(1-t),\varphi_{t},\psi_{t})\\
\end{split}.
\end{equation}
Now, we define:
\[\Phi_{1}^{t}:H_{1}\longrightarrow \mathbb{R}^{2}\times (S^{1})^{2}\]
by:
\begin{equation}\label{eq:t2}
\Phi_{1}^{t}: P\longmapsto
\begin{cases}
\Phi_{1T}^{t}(P)&\text{ if } P\in H_{1T};\\&\\
\Phi_{1L}^{t}(P) &\text{ if } P\in H_{1L}.
\end{cases}
\end{equation}

\end{defn}

\begin{claim}
The map $\Phi_{1}^{t}$ is well defined.
\end{claim}
\proof
On the subset $S:=H_{1T}\cap H_{1L}=\left\{(x,y,\varphi,\psi)\in H_{1}|y=2x+\ln2\right\}$, the two equations defining $\Phi_{1}^{t}$ agree.
Indeed, they clearly agree on the first two components, for all $(x,y,\varphi,\psi)\in S$. They also agree on the last two components as, by Claim~$\ref{cl:7}$, the definitions of $a$ and $a'$ agree on $S$.
\endproof

\begin{prop}
Let $i_{H_{1}}:H_{1}\hookrightarrow(\mathbb{C}^{*})^{2}$ be the canonical inclusion.
We have $\Phi_{1}^{0}=i_{H_{1}}$ and $\Phi_{1}^{1}(H_{1})=H_{1trop}$.
\end{prop}
\proof
The first part of the statement follows immediately by substituting $t=~0$ in the definition of $\Phi^{t}_{1}$. From the definition of $\Phi_{1T}^{1}$  and from Lemma $\ref{lem:8}$ we get $\pi(\Phi_{1T}^{1}(H_{1T}))=(0,0)$ and $Arg(\Phi_{1T}^{1}(H_{1T}))=(\varphi_{1},\psi_{1})(H_{1T})=\overline{Arg(H_{1})}$, thus\\ $\Phi_{1T}^{1}(H_{1T})= H_{1tropT}$. Now, let $c\geq\ln2$. Consider the subset of $H_{1L}$:
\[W_{c}:=\left\{(x,y,\varphi,\psi)\in H_{1L}|y=2x+c\right\}.\]
From the definition of $\Phi_{1L}^{1}$ we get:
\[\Phi_{1L}^{1}(W_{c})=\left\{(x-\frac{2x+c-\ln2}{2},0,\varphi,\pi)|\varphi\in S^{1}\right\}\]
that is:
\[\Phi_{1L}^{1}(W_{c})=\left\{(\frac{-c+\ln2}{2},0,\varphi,\pi)|\varphi\in S^{1}\right\}.\]
Since we have:
\[H_{1L}=\bigsqcup_{c\geq\ln2}W_{c},\]
we get:
\[\Phi_{1L}^{1}(H_{1L})=\left\{(\frac{-c+\ln2}{2},0,\varphi,\pi)|\varphi\in S^{1}, c\geq\ln2\right\}.\]
From $c\geq\ln2$, we get $\frac{-c+\ln2}{2}\leq0$. Therefore, $\Phi_{1L}^{1}(H_{1L})=H_{1tropL}$.

\endproof

\begin{prop}\label{prop:34}
The map $\Phi_{1}^{t}$ is injective.
\end{prop}
\proof

Let $P,P'\in H_{1T}$ and assume $\Phi_{1T}^{t}(P)=\Phi_{1T}^{t}(P')$. By Lemma $\ref{lem:8}$, the map $(\varphi_{t},\psi_{t})$ is injective on $Arg(H_{1T})$, hence $\varphi=\varphi'$ and $\psi=\psi'$ and using Lemma $\ref{lem:9}$ we get $P=P'$. Now, let $P,P'\in H_{1L}$ and assume $\Phi_{1L}^{t}(P)=\Phi_{1L}^{t}(P')$. If $t\neq 1$, then $x=x'$ and $y=y'$. Therefore, by Proposition $\ref{prop:2}$ $(\varphi,\psi)=(\varphi',\psi')$ or $(\varphi,\psi)= (2\pi-\varphi',2\pi-\psi')$. If both hold, then we get $(\varphi,\psi)=(\varphi',\psi')=(\pi,\pi)$ and hence $P=P'$. It cannot happen that only the second equality holds as by definition of $(\varphi_{t},\psi_{t})$, $\Phi_{1L}^{t}(P)=\Phi_{1L}^{t}(P')$ implies that $(\varphi,\psi)$ and $(\varphi',\psi')$ are either both in $A_{1}$ or both in $A_{2}$. If $t=1$, then $\Phi_{1L}^{1}(P)=\Phi_{1L}^{1}(P')$ implies that $(\varphi,\psi)$ and $(\varphi',\psi')$ lie both on the same straight line passing through $O$ or through $O'$. On the other hand, we also get $2(x-x')=y-y'$, which means that $P,P'\in W_{c}$,  where $c\geq ln 2$. Thus, by Proposition $\ref{prop:4}$ we get $(\varphi,\psi)=(\varphi',\psi')$, and hence $P=P'$ as $Arg|_{W_{c}}$ is injective.
\endproof
We recall the following well-known results in general topology:
\begin{lem}\label{lem:10}
Let $f:X\longrightarrow Y$ be a continuous map between locally compact Hausdorff spaces. If $f$ is proper, then it is closed.
\end{lem}
\begin{lem}\label{lem:11}
Let $f:X\longrightarrow Y$ be a bijective continuous map. If $f$ is closed, then it is bicontinuous.
\end{lem}

\begin{lem}\label{lem:12}
Let $p:X\times Y\longrightarrow X$ be the projection onto the first factor. If $Y$ is compact, then $p$ is proper.
\end{lem}

Now, we prove:
\begin{lem}\label{lem:13}
For any $t\in\left[0,1\right]$, the map $\Phi_{1L}^{t}:H_{1L}\rightarrow \Phi_{1L}^{t}(H_{1L})$ is a proper map between locally compact Hausdorff spaces.
\end{lem}
\proof
For any $t\in\left[0,1\right]$, $\Phi_{1L}^{t}(H_{1L})$ is clearly a closed subset of $\mathbb{R}^{2}\times (S^{1})^{2}$, hence it is locally compact and Husdorff. Now, we show that 
$\Phi_{1L}^{t}$ is proper.\\Let $\pi:\mathbb{R}^{2}\times (S^{1})^{2}\rightarrow \mathbb{R}^{2}$ be the projection and consider the family of maps:
\begin{alignat*}{9}
g_{t}:&&\mathbb{R}^{2}&\longrightarrow\mathbb{R}^{2} \\
      &&(x,y)         &\longmapsto(x-t\frac{y-\ln2}{2},y(1-t)),
\end{alignat*}
for $t\in\left[0,1\right]$. We have the following commutative diagram:

\[	\begin{tikzpicture}[descr/.style={fill=white,inner 		
					sep=2.5pt}]
					\matrix (m) [matrix of math nodes, row sep=3em,
					column sep=3em]
					{ H_{1L} &  \Phi_{1L}^{t}(H_{1L})\\
					   \pi(H_{1L}) & \pi(\Phi_{1L}^{t}(H_{1L})) \\ };
					\path[->,font=\scriptsize]
					(m-1-1) edge node[above] {$\Phi_{1L}^{t}$} (m-1-2) 
					(m-1-1) edge node[left] {$\pi|_{H_{1L}}$} (m-2-1)
					(m-1-2) edge node[right] {$\pi|_{\Phi_{1L}^{t}(H_{1L})}$} (m-2-2)
					(m-2-1) edge node[above] {$g_{t}|_{\pi(H_{1L})}$} (m-2-2);
				\end{tikzpicture}.\]
For $t\neq 1$, the map $g_{t}$ is a homeomorphism, hence $g_{t}|_{\pi(H_{1L})}$ is proper. If $t=1$, then we see that $g_{1}$ projects the points on the straight line $y=2x+c$, for $c\in\mathbb{R}$, to the point $(\frac{-c+\ln2}{2},0)$. Thus, $g_{1}|_{\pi(H_{1L})}$ is also proper. Now, by Lemma $\ref{lem:12}$, $\pi$ is a proper map and since $\Phi_{1L}^{t}(H_{1L})$ is a closed subset of $\mathbb{R}^{2}\times (S^{1})^{2}$, we have that $\pi|_{\Phi^{1L}_{t}(H_{1L})}$ is still a proper map for all $t\in\left[0,1\right]$. Thus, $g_{t}|_{\pi(H_{1L})}\circ\pi|_{H_{1L}}$ is a proper map. But, $g_{t}|_{\pi(H_{1L})}\circ\pi|_{H_{1L}}=\pi|_{\Phi_{1L}^{t}(H_{1L})}\circ\Phi_{1L}^{t}$, one easily sees that $\Phi_{1L}^{t}$ is a proper map too.
\endproof

\begin{thm}\label{thm:1}
The map $\Phi_{1}^{t}:H_{1}\longrightarrow \Phi_{1}^{t}(H_{1})$ is a homeomorphism.
\end{thm}
\proof
The map is bijective by Proposition $\ref{prop:34}$ and it is clearly continuous. The map $\Phi_{1T}^{t}:H_{1T}\rightarrow \Phi_{1T}^{t}(H_{1T})$ is bicontinuous as it is a bijective continuous map from a compact space onto a Hausdorff space.
We also have that $\Phi_{1L}^{t}:H_{1L}\rightarrow \Phi_{1L}^{t}(H_{1L})$ is bicontinuous. Indeed, by Lemma $\ref{lem:13}$, $\Phi_{1L}^{t}$ is a proper map between locally compact Hausdorff spaces. Therefore, by Lemma $\ref{lem:10}$, $\Phi_{1L}^{t}$ is a closed map and hence, by Lemma $\ref{lem:11}$, it is bicontinuous.
\endproof

\begin{thm}
Let $H\subset(\mathbb{C}^{*})^{2}$ the complex line $1+z_{1}+z_{2}=0$ and $H_{trop}\subset(\mathbb{C}^{*})^{2}$ its associated phase tropical line. Let $i_{H}:H\hookrightarrow(\mathbb{C}^{*})^{2}$ be the canonical embedding.
Then there exists a continuous map: 
\[\Psi:H\times [0,1]\longrightarrow (\mathbb{C}^{*})^{2}\] 
such that the family of maps:  
\[\Psi_{t}:H\longrightarrow (\mathbb{C}^{*})^{2}, P\longmapsto \Psi(P,t),\]
for $t\in [0,1]$, has the following properties: 
\begin{enumerate}
	\item $\Psi_{0}=i_{H}$;
  \item $\Psi_{1}(H)=H_{trop}$;
	\item $\Psi_{t}$ is a homeomorphism onto the image,  for each $t\in\left[0,1\right]$.
\end{enumerate}

\end{thm} 

\proof

Let $t\in\left[0,1\right]$, $\Phi^{t}_{1}$ as in Definition $\ref{defn:6}$ and Theorem $\ref{thm:1}$, let $\lambda$ be as in Lemma $\ref{lem:4}$ and $H_{2},H_{3}\subset H$ as in Definition $\ref{defn:4}$. We can define two maps:
\[\Phi^{t}_{2}:H_{2}\longrightarrow \mathbb{R}^{2}\times (S^{1})^{2}, \Phi^{t}_{2}:=\lambda\circ\Phi_{1}^{t}\circ\lambda^{-1}\]
and
\[\Phi^{t}_{3}:H_{3}\longrightarrow \mathbb{R}^{2}\times (S^{1})^{2}, \Phi^{t}_{3}:=\lambda^{-1}\circ\Phi_{1}^{t}\circ\lambda.\]

Using Lemma $\ref{lem:4}$ we easily see that $\Phi^{0}_{2}=i_{H_{2}}$, $\Phi^{0}_{3}=i_{H_{3}}$ and $\Phi^{1}_{2}(H_{2})=H_{2trop}$, $\Phi^{1}_{3}(H_{3})=H_{3trop}$. Moreover, $\Phi^{t}_{2}$ and $\Phi^{t}_{3}$ are homeomorphisms onto their images being compositions of homeomorphisms. We define a map:
\[\Psi'_{t}:H\longrightarrow \mathbb{R}^{2}\times(S^{1})^{2}\]
via
\begin{equation}
\Psi'_{t}:P\longmapsto
\begin{cases}
\Phi^{t}_{1}(P),&\text{ if } P\in H_{1}\\
\Phi^{t}_{2}(P),&\text{ if } P\in H_{2}\\
\Phi^{t}_{3}(P),&\text{ if } P\in H_{3}
\end{cases}.
\end{equation}
It is clearly well defined. Let us consider the homeomorphism:
\begin{equation}
\begin{split}
h:(\mathbb{C}^{*})^{2}&\longrightarrow \mathbb{R}^{2}\times (S^{1})^{2}\\
(z_1,z_2)&\longmapsto(\ln|z_1|,\ln|z_2|,\arg z_1,\arg z_2)\\
\end{split}.
\end{equation}
The map
\[\Psi_{t}:= h^{-1}\circ\Psi'_{t}\circ h\]
satisfies properties \textit{1.-3.} in the statement of the theorem. Now, 
\[\Psi:H\times\left[0,1\right]\longrightarrow (\mathbb{C}^{*})^{2}, (P,t)\longmapsto \Psi_{t}(P)\]
is the claimed map.

\endproof

\cleardoublepage
\addcontentsline{toc}{chapter}{References}	
\bibliographystyle{alpha}
\bibliography{bibliography2}

\end{document}